\documentclass{amsart}
\usepackage[centertags]{amsmath}
\usepackage{amsfonts}
\usepackage{amsthm}
\usepackage{newlfont}
\usepackage{amscd}
\usepackage{amsgen}
\usepackage{pb-diagram}
\newlength{\defbaselineskip} 
\theoremstyle{plain}
\newtheorem{thm}{Theorem}[section]
\newtheorem{prop}[thm]{Proposition}
\newtheorem{cor}[thm]{Corollary}

\newtheorem{lem}[thm]{Lemma}

\theoremstyle{remark} \newtheorem{rem}[thm]{Remark}
\theoremstyle{definition} 
\theoremstyle{definition} \newtheorem{ex}[thm]{Example} %

 \numberwithin{equation}{section}
\newcommand{\op}{\operatorname}

\newcommand{\ra}{\rightarrow}

\begin{document}
\title{Fiber products of elliptic surfaces with section and associated Kummer fibrations.}
\author{Grzegorz Kapustka and Micha{\l}\ Kapustka}
\thanks{Mathematical subject classification 14J32; 14J27}
\thanks{The project is co-financed by the EU funds and national budjet.}
\maketitle \begin{abstract}{ We investigate Calabi--Yau three
folds which are small resolutions of fiber products of elliptic
surfaces with section admitting reduced fibers. We start by the
classification of all fibers that can appear on such varieties.
Then, we find formulas to compute the~Hodge numbers of obtained
three folds in terms of the types of singular fibers of
the~elli\-ptic surfaces. Next we study Kummer fibrations
associated to these fiber products.}
\end{abstract}
\section*{Introduction}

Calabi--Yau three folds which are small resolutions of nodal fiber
products of relatively minimal rational elliptic surfaces with
section were introduced by Schoen in \cite{Schoen}. This class of
manifolds appeared to be a very good background for~studying for
instance arithmetics and modularity. For instance, most of the
examples of non rigid modular Calabi--Yau three folds are
connected in some way with this family.

The main goal of this article is to extend
the family studied by Schoen.

In the first part of the paper we investigate Calabi--Yau three
folds which are small resolutions of fiber products of elliptic
surfaces with section admitting reduced fibers. We allow all types
and configurations of such fibers on the rational elliptic
surfaces and hence obtain more complicated rational double points
to resolve. We start by the classification of all fibers that can
appear on such varieties. Then, we find formulas to compute
the~Hodge numbers of obtained three folds in terms of the types of
singular fibers of the~elli\-ptic surfaces. Next we focus on
finding geometric interpretation of the deformations of studied
manifolds. In the general context (when the fibers of the elliptic
surfaces can be nonreduced) we prove that a generic deformation of
a Calabi--Yau resolution of a fiber product of two elliptic
surfaces with section is birational to such a fiber product. In
more special cases we give two ways of interpretation. The first
is in terms of deformations of the configurations of singular
fibers on both elliptic surfaces. The second bases on the
description of studied three folds as iterated covers of
$\mathbb{P}^3$ branched along two quartic cones.

The rest of the paper is devoted to the study of Kummer
fibrations associated to these fiber products. Their construction
consists of taking the quotient of a fiber product of two elliptic
surfaces with section by a suitable involution and next resolving
obtained singularities. The described three folds admit fibrations by Kummer
surfaces. In general fibred manifolds play a special role in arithmetics. Moreover, general three folds
fibred by K3 surfaces are still not well understood. The advantage of
the Kummer fibrations constructed is that they belong to the well
described family of double octics. This enables us to study the Hodge numbers of these
varieties.

\textbf{Acknowledgements} We would like to thank our advisor Dr hab. S. Cynk for his constant
enormous help during the work on this paper. We are grateful to Prof. A. Langer for his precious remarks. We would like also to thank the Referee of this 
paper for many interesting suggestions.

\section{The fiber products of elliptic surfaces with section} In
\cite{Schoen} Schoen introduced a new class of Calabi--Yau
varieties. These are small resolutions of nodal fiber products of
relatively minimal elliptic surfaces with section. This family was
investigated by many algebraic geometers. In this section we would
like to extend Schoen's construction and study obtained varieties
a bit closer.

In this paper by the term elliptic surface we mean a surface $S$
together with a morphism $\pi: S\ra\mathbb{P}^1$ such that the
generic fiber of $\pi$ is a smooth elliptic curve. By abuse of
notation we denote such elliptic surface simply by $S$. Using this
term we moreover require that $S$ is minimal in the sense that it
does not have -1 curves as components of fibers. An elliptic
surface $S$ admits a section if the associated morphism $\pi$
admits a section. For simplicity we identify the section with its
image.

Let $S_1$ and $S_2$ be two rational elliptic surfaces with section
and assume that they have only reduced fibers. Let $X$ be the
fiber product $S_1 \times_{\mathbb{P}^1} S_2$. If we see it as the
inverse image of the diagonal in $\mathbb{P}^1 \times
\mathbb{P}^1$ to the product $S_1 \times S_2$, using the
adjunction formula (for more details see \cite{Schoen}) we prove
that $X$ is a singular Calabi--Yau three fold. Hence, if $X$
admits a small resolution $\hat{X}$, then $\hat{X}$ is a
Calabi--Yau three fold. In this section we are interested in
studying Hodge numbers of such Calabi--Yau three folds.

Let us start with the following proposition.
\begin{prop} \label{lem male rozwiazania}Let $S_1$ and $S_2$ be two rational
elliptic surfaces with section. Then the three fold $X=S_1
\times_{\mathbb{P}^1} S_2$ admits a small resolution if and only
if all fibers of $X$ are of one of the following types: $F\times
I_0$, $I_n \times I_m$, $III \times I_n$, $III \times III$, $IV
\times I_n$, $II\times II$, where $F$ denotes any fiber of an
elliptic ruled surfaces.
\end{prop}

\begin{proof}
Observe that all the singularities of $X$ are lying on a fiber which
is the product of two singular fibers. We can find explicitly the
local equations of the variety around each of these singularities
and check which ones admit a small resolution.

The local equations of the surfaces around a singularity of the
fiber are the following:
\begin{itemize}
\item[$I_n$] The local equation around each singular point of the
fiber is $t=xy$.
 \item[$II$] The local equation around the singularity of the fiber is
 $t=y^2-x^3$.
 \item[$III$] The local equation around the singularity of the fiber is $t=
 x(y^2-x)$.
 \item[$IV$] The local equation around the singularity of the fiber is
 $t=xy(x+y)$.
\end{itemize}
In each of these cases the fibration is given by the projection
onto the coordinate $t$. It follows that the fiber product around
each of its singularities has one of the following local
equations.
\begin{itemize}
\item[$I_n \times I_m$] The local equation around each of the
$n\cdot m$ singularities is $xy=u^2+v^2$. This is an equation
defining a node, which admits a small resolution.
 \item[$I_n \times II$] The local equation around each of the $n$
 singularities is $xy=u^2-v^3$. This is an equation
 of an $A_2$ singularity. It does not admit a small
 resolution as it is factorial.
 \item[$I_n \times III $] The local equation around each of the $n$
 singularities is $xy=u(v^2-u)$. This equation defines
 an $A_3$ singularity. It admits a small
 resolution which is given by the blowing up of the plane $\{x=u=0\}$.
 \item[$I_n \times IV $] The local equation around each of the $n$
 singularities is $xy=uv(u+v)$. This is a
 singularity of type $D_4$. It admit a small
 resolution which is given by the consecutive blowing up of the plane $\{x=u=0\}$
 and the proper transform of $\{y=v=0\}$.
 \item[$II \times II$] The local equation around the
 singularity is $x^2-y^3=u^2-v^3$. This is also a
 singularity $D_4$. Its small resolution is given by the blow up
 of the plane $x-u=y-v=0$ and then taking a small resolution of the resulting node.
 \item[$II \times III$] The local equation around the
 singularity is $x^2-y^3=u(v^2-u)$. This is a
 singularity of type $E_6$, which is factorial from \cite[thm B]{hui lin}. Hence it does not admit any small resolution.

 \item[$II \times IV$] The local equation around the
 singularity is $x^2-y^3=uv(u+v)$. It is a double cover of $\mathbb{C}^3$
 ramified over a triple point. From \cite[thm B]{hui lin} it is a factorial
 singularity. It does not admit any small resolution.
 \item[$III \times III $] The local equation around the
 singularity is $x(y^2-x)=u(v^2-u)$. It admits a small
 resolution which is given by the blowing up of the plane $x=u=0$ and a small resolution of resulting nodes.
 \item[$III \times IV $] The local equation around the
 singularity is $x(y^2-x)=uv(u+v)$. This is a
 singularity that does not admit any small resolution. We can see this
 in the following way. After blowing up the plane $x=u=0$ we
 obtain an $A_2$ hence factorial singularity. Now we can use an
 argument of Koll\'ar to prove that the small resolution does not
 exist. Indeed, if we had such a resolution, then it would differ
 from the obtained $A_2$ by a sequence of flops \cite[theorem 4.9]{Kollarflop}. This gives a
 contradiction with the fact that flops in dimension 3 do not change the types of singularities \cite[theorem 2.4]{Kollarflop}.
 \item[$IV \times IV $] The local equation around the
 singularity is of type $xy(x+y)=uv(u+v)$. This is a simple triple point. It admits
 a big crepant resolution, hence do not have any small one.
\end{itemize}
\end{proof}

\begin{prop}\label{projective fiber product} The three fold $X$ admits a projective small resolution $\hat{X}$ if
moreover no fiber is of type $II\times II$, $III\times III$ or
$I_1\times F$, where $F$ is a singular fiber.
\end{prop}
\begin{proof}
We easily check that for each of allowed singularities the
analytical small resolution described above is given by a sequence
of blowings up of components of fibers (they are all isomorphic to
$\mathbb{P}^1 \times \mathbb{P}^1$).
\end{proof}
To study Hodge numbers of constructed varieties we introduce some
notation. Let $A_1$,$A_2\subset \mathbb{P}^1$ be the sets of
singular fibers of surfaces $S_1$ and $S_2$ respectively. Let
$S''=A_1 \cap A_2$. Let $b_t$ and $b'_t$ for $t\in \mathbb{P}^1$
denote the number of components of the fibers over the point $t$
of $S_1$ and $S_2$ respectively.

Let us start by computing $h^{1,1}(\hat{X})$. We adapt the
arguments of \cite{Schoen} to our wider context.  We have the
exact sequence
$$\text{\begin{minipage}{5.5cm}{\begin{center}The free abelian group generated by all components
of the fibers \end{center}}
\end{minipage}}\longrightarrow\operatorname{Pic}(X)\longrightarrow
\operatorname{Pic}(X_\eta)\longrightarrow 0$$
 where $X_{\eta}$ denotes the generic fiber of $X$ in the sense of
 schemes. Moreover, all relations between images of components of the fibers are
 induced by the linear equivalence of the fibers.
We get

$$\op
{rk}\operatorname{Pic}(\hat{X})=1+\sum_{t \in A_1\cup A_2} (b_t b'_t
-1)+\op{rk}\operatorname{Pic}(\hat{X}_{\eta}).$$

Continuing, define $d(X)=1$ if and only if $S_1$ and $S_2$ are
isogenous and $d(X)=0$ otherwise. Let $S_{1\eta}$, $S_{2\eta}$
denote the generic (in the sense of schemes) fibers of the
respective surfaces $S_1$ and $S_2$
$$\op{rk}\operatorname{Pic}(\hat{X}_{\eta})=\op{rk}\operatorname{Pic} (S_{1\eta}) + \op{rk}\operatorname{Pic} (S_{2\eta})+d(X)=
18-\sum_{t\in A_1} (b_t-1) -\sum_{t\in A_2} (b'_t-1)+d(X).$$

Together this gives
$$ h^{1,1}(\hat{X})=\sum_{t\in A_1\cup A_2} (b_t b'_t -1)+ 19-\sum_{t\in A_1} (b_t-1) -\sum_{t\in A_2} (b'_t-1)+d(X).$$

Now, to get a formula for $h^{1,2}(\hat{X})$ we compute the Euler
characteristic of $\hat{X}$. First for $t\in\mathbb{P}^1$ define
$F_t$ and $F'_t$ to be the fibers over $t$ of surfaces $S_1$ and
$S_2$ respectively. For all $s \in \operatorname{Sing}(X)$ define
$E_s$ to be the exceptional curve of the small resolution of the
singularity in $s$. Then
$$\chi(\hat{X})= \chi(X) + \sum_{s\in \op{Sing}(X)} (\chi(E_s)-1) $$
and
$$\chi(X)=\sum_{t \in S''} \chi(F_t) \chi(F'_t). $$

Hence, to compute the Euler characteristic $\chi(\hat{X})$ we need
only to compute the Euler characteristics of the fibers of the
surfaces, and of the exceptional curves of each resolved
singularity. We already have
$$\chi(F)=
\begin{cases}
 n \text{ for } F \text{ of type } I_n\\
 3 \text{ for } F \text{ of type } III\\
 4 \text{ for } F \text{ of type } IV.
\end{cases}
$$
To compute the Euler characteristic of an exceptional curve we
only need to know the number of local blowings up of Weil divisors
we performed during the resolution of singularities described in
the proof of Proposition \ref{lem male rozwiazania}. This gives
$$\chi(E_i)=
\begin{cases}
 2 \text{ for all nodes and the singularities of $I_n \times III$}\\
 3 \text{ for the singularities of $IV \times I_n$}
\end{cases}
$$
Together this simplifies to a formula for the Euler characteristic
depending only on the numbers of components of fibers
$$\chi(\hat{X})=2\sum_{t \in S''} b_t b'_t.$$
To get the formula for $h^{1,2}$ we use
$$h^{1,2}(\hat{X})= h^{1,1}(\hat{X})-\frac {\chi(\hat{X})}{2}.$$
Finally we obtain
$$\begin{aligned}
h^{1,2}(\hat{X})= 19+d(X)-\sum_{t \in S''} (b_t + b'_t -1).
\end{aligned}
 $$
\begin{rem} We have computed Hodge numbers for the class of projective
Calabi--Yau three folds described in Proposition \ref{projective
fiber product}. Let us remind that any, including the
non-projective, introduced Calabi--Yau small resolution is a
Moishezon manifold. We can hence study its Hodge numbers.
Moreover, in our case we can extend the above computations to all
Calabi--Yau three folds $\hat{X}$. The only difference is that we
need to put appropriate Euler characteristics of the exceptional
divisors corresponding to obtained singularities on $X$. This is
not hard as these depend only on the number of local blowings up
in the resolutions described by the proof of Proposition \ref{lem
male rozwiazania}.
\end{rem}
On any smooth Calabi--Yau three fold $X$ the number $h^{1,2}(X)$
is equal to the space of infinitesimal and by the theorem of
Bogomolov--Tian--Todorov also small deformations of this three
fold. Let us focus on this interpretation and find a natural local
description of the Kuranishi space of $\hat{X}$.

We shall first present some argument concerning deformations of
general Calabi--Yau three folds which are resolutions of fiber
products of elliptic surfaces with section.
\subsection{Deformations in the general case}
In this subsection we do not restrict the type of fibers of the
product and prove the following
\begin{thm}\label{infinitesimal deformations of fiber products}
Let $X$ be a Calabi--Yau three fold, which is a resolution of a
fiber product of two rational elliptic surfaces with section. Then
a generic deformation of $X$ is birational to a fiber product of
this type.
\end{thm}
The proof is based on the following lemma.
\begin{lem} \label{abel}Let $X$ be a Calabi--Yau manifold containing the
product of two elliptic curves $E_1$ and $E_2$. Assume moreover
that $h^0(E_1,N_{E_1|X})=h^0(E_2,N_{E_2|X})=2$. Then $X$ is
birational to a fiber product of rational elliptic surfaces with
section.
\end{lem}

\begin{proof}[Proof of Lemma \ref{abel}]
Let $A= E_1 \times E_2$, where $E_1$ and $E_2$ are elliptic
curves, be the abelian surface contained in $X$. We can easily
compute that the linear system $|A|$ is one dimensional and next
we obtain $X$ to admit an abelian fibration over $\mathbb{P}^1$
with special fiber $A$ (For more details see \cite{GP}).

The infinitesimal deformations of $A$ in $X$ are induced by the
morphism
$$ H^0(A,N_{A|X})\longrightarrow H^1(A,T_A)$$
coming from the exact sequence
\begin{equation}\label{exact sequenc TA}
0\longrightarrow T_{A} \longrightarrow T_{X}|_A \longrightarrow
N_{A|X}\longrightarrow 0.
\end{equation}
We will prove that these deformations of $A$ in $X$ induce
deformations of pairs $(A,E_1)$ and $(A,E_2)$. Consequently a
generic abelian fiber of the fibration will contain two elliptic
curves intersecting in one point, thus will be a product.

Let us concentrate on the deformation of $(A,E_1)$. We have two
exact sequences
$$0\longrightarrow T_{E_1} \longrightarrow T_{A}|_{E_1}
\longrightarrow N_{E_1|A}\longrightarrow 0.
$$
$$0\longrightarrow T_{A}(-E_1) \longrightarrow T_{A}
\longrightarrow T_{A}|_{E_1} \longrightarrow 0.
$$
Together with the exact sequence \ref{exact sequenc TA} they
induce a diagram in cohomology.
$$
\begin{CD}  H^0(A,N_{A|X})@>>> H^1(A,T_A)@.\\
@. @ VVV @.\\
 H^1(E_1,T_{E_1})@ >>> H^1(E_1,T_A|_{E_1})@>>>H^1(E_1,N_{E_1|A})\\
\end{CD} $$

We need to prove that there are no obstructions to lift the
considered infinitesimal deformations of $A$ in $X$ to
deformations of $E_1$ in $X$, i.e. that the composition morphism
$$ H^0(A,N_{A|X})\longrightarrow H^1(A,T_A)\longrightarrow H^1(E_1,T_A|_{E_1})\longrightarrow
H^1(E_1,N_{E_1|A}),$$  is trivial. Observe that this morphism
factorizes by the suitable morphism induced by the following exact
sequence:
\begin{equation}
\label{exact sequenc NE1} 0\longrightarrow N_{E_1 | A}
\longrightarrow N_{E_1 |X} \longrightarrow
N_{A|X}|_{E_1}\longrightarrow 0
\end{equation}
From the adjunction formula $N_{E_1 |A} \simeq N_{A |X}|_{E_1}
\simeq \mathcal{O}_{E_1}$. Now from the assumption and the
associated cohomology sequence
$$0\longrightarrow H^0(E_1,N_{E_1 |A}) \longrightarrow H^0(E_1,N_{E_1 |X})\longrightarrow H^0(E_1,N_{A
|X}|_{E_1})\longrightarrow H^1(E_1, N_{E_1|A})$$ for dimension
reasons the morphism $H^0(E_1,N_{E_1 |X})\longrightarrow
H^0(E_1,N_{A |X}|_{E_1})$ is surjective. Hence $H^0(E_1,N_{A
|X}|_{E_1})\longrightarrow H^1(E_1, N_{E_1|A})$ is trivial.

To finish the proof we remind that deformations of abelian
surfaces on $X$ as well as deformations of elliptic curves are
unobstructed. Hence the obtained deformations of elliptic curves
$E_1$ and $E_2$ sweep out two elliptic surfaces $S_1$ and $S_2$.
The three fold $X$ is then birational to the fiber product $S_1
\times_{\mathbb{P}_1}S_2$. The two intersecting elliptic curves in
each smooth fiber yield sections in both elliptic surfaces.
Moreover, as the deformation of a Calabi--Yau three fold is a
Calabi--Yau three fold the two surfaces $S_1$ and $S_2$ are
rational (see \cite{Schoen}).

\end{proof}
We are now ready to prove the main theorem:
\begin{proof}[Proof of Theorem \ref{infinitesimal deformations of fiber products}]
 Let $X$ be a resolution of singularities of a
 fiber product $Y$ of two rational elliptic surfaces with section.
Then $X$ contains a surface $A$ which is a product of two elliptic
curves. By \cite{Wilson2} the K\"ahler cone of $X$ and the
K\"ahler cone of a generic deformation $X'$ of $X$ can be seen as
subcones of the same space $H^2(X,\mathbb{R})$. From the same
paper we know that each abelian fibration is given by a divisor
$D_A$ such that $D_A^3=0$, $D_A^2\equiv 0$, $c_2(X).D_A =0$ and
each elliptic fibration is given by $D_E$ such that $D_A^3=0$,
$D_A^2 \not\equiv 0$, $c_2(X).D_A \geq 0$. Moreover by
\cite{Wilson4} the Kahler cone of $X$ is a subcone of the Kahler
cone of $X'$ cut out by some half space. By \cite[cor. 2, rem.
1]{Wilson5} the divisors $D_{E_1}$ and $D_{E_2}$ induce also two
elliptic fibrations on $X'$.

We are interested in finding an abelian surface containing two
elliptic curves each from a different fibration. Let $S$ be a
section of the fibration given by $D_{E_1}$ on $X$ (such section
exists as $X$ is a fiber product of two elliptic surfaces with
section). Then $S$ is a rational elliptic surface. In particular
$H^1(S,\mathcal{O}_S)=0$, hence any deformation $X'$ of $X$
induces a deformation $S'$ of $S$. As the intersection form is
fixed in the deformation, $S'$ is covered by fibers of the
fibration given by $D_{E_1}$ and cuts the fibers of the fibration
given by $D_{E_2}$ in one point. Thus we can consider the surface
$A'$ on $X'$ swept out by curves from one fibration over a curve
from the second. Then $A'^2\equiv 0$, $A^3=0$ and $A'$ is not a K3
surface (it admits an elliptic fibration over an elliptic curve).
This implies that $A$ is an abelian surface on $X'$ containing two
elliptic curves $E'_1$, $E'_2$ meeting in one point. Hence it is
the product of these two curves. Moreover as both curves come from
two-dimensional fibrations we have
$h^0(E'_1,N_{E'_1|X'})=h^0(E'_2,N_{E'_2|X'})\geq 2$. By lemma
\ref{abel} this ends the proof.
\end{proof}

Let us moreover point out the following
\begin{rem} If a ruled surface over a curve of genus $g\geq 1$ is
contained in $X$, then it has to be a component of a fiber
 of type $I_0 \times F$.

 Indeed suppose to the contrary. Let $S$ be such a surface contained in
 $X$. Then there exists a surjection from $S$ to one of the elliptic
 surfaces (as the fiber product of two multi-sections is a multi-section).
 The inverse-image $F$ of the generic elliptic fiber have to be an effective
 smooth (as it is base point free) divisor with self-intersection 0. Hence one of the following holds:
\begin{itemize}
\item The divisor $F$ is the sum of some fibers of
 the ruled surface. This is impossible as there exists no morphism
 from $\mathbb{P}^1$ to a curve of positive genus (L\"uroth's theorem).
\item The surface $S$ is a product of a genus $g$ curve with a
line. This is impossible since the inverse image of the section of
the elliptic surface would have to be an exceptional curve. Such
curves do not exist on products.
\end{itemize}
\end{rem}

\begin{cor}\label{generic F not I1,II } A generic Calabi--Yau manifold from the Kuranishi
space of a fiber product of two rational elliptic surfaces with
section is also a fiber product of two rational elliptic surfaces
with section and has no fibers of type
 $I_0 \times F$ for $F \not= I_1, II$
\end{cor}

\subsection{Deformations of constructed three folds}
In this subsection we study projective Calabi--Yau three folds $X$
introduced in Proposition \ref{projective fiber product}. We are
looking for a more precise description of the space of small
deformations of $\hat{X}$. That is we would like to see
deformations from the point of view of singular fibers. Let us
first find a generic variety obtained as a deformation of
$\hat{X}$.

\begin{prop} \label{o postaci deformacji produktu}
A generic small deformation $\hat{Y}$ of $\hat{X}$ is a small
resolution of a fiber product $Y$ of two elliptic surfaces $R_1$
and $R_2$ which are small deformations of $S_1$ and $S_2$
respectively and such that each fiber of $Y$ corresponds to a
fiber of $X$ with the same number of components. Moreover, such
variety $Y$ has no fibers of type $F \times I_0$ for $F\not=I_1
\text{ or } II$.
\end{prop}
\begin{proof} Let us consider the Weierstrass type equations of the
surfaces $S_1=\{y^2=x^3+f(t)x+g(t)\}$ and
$S_2=\{y'^2=x'^3+f'(t)x'+g'(t)\}$, where $f$, $g$, $f'$, $g'$ are
polynomials of degrees 4,6,4,6 respectively. Consider the set $Q$
of all quadruples of polynomials $F$, $G$, $F'$, $G'$, of degrees
4,6,4,6 respectively, satisfying the following conditions
\begin{enumerate}
 \item The pairs $F$, $G$ and $F'$, $G'$ define by the
Weierstrass type equations two elliptic surfaces $R_1$ and $R_2$.
 \item There is a bijection $\beta$ between the set $R''$ of common
singular fibers of $R_1$ and $R_2$ and the corresponding set
$S''$.
 \item for all $t\in R''$ the number of components of
the fibers over $t$ of both surfaces $R_1$ and $R_2$ are the same
as the numbers of components of fibers of respective surfaces
$S_1$ and $S_2$ over $\beta(t)$.
\end{enumerate}
\textbf{Claim} The set $Q$ is an algebraic set of dimension at
least $h^{1,2}(\hat{X})+5$.

The set is clearly nonempty. Hence to compute its dimension we
compute the number of conditions we impose on the set of all
quadruples of polynomials of respective degrees 4,6,4,6. Note that
requiring that a pair $(F,G)$ defines a surface with a fiber with
$b_t$ components over a chosen $t$ is imposing $b_t$ conditions.
Now:
\begin{itemize}
 \item In the case where the surfaces $S_1$ and $S_2$ are not
isogenous the number of conditions is $\sum_{t\in R''}
(b_t+b'_t-1)$.
 \item In the case where the surfaces $S_1$ and $S_2$
are isogenous the number of conditions describing isogenous pairs
of surfaces is $\sum_{t\in R''} (b_t+b'_t-1)-1$.
\end{itemize}
The claim is proved.

Element from the set $Q$ correspond to elliptic surfaces with
fibers having the same number of components over $\beta(t)$ as the
original surface in $t$. Let us consider the component of $Q$
containing the quadruple describing $X$. Using the fact that an
elliptic rational surface has only one Weierstrass representation
up to proportion and after dividing by the automorphism group of
$\mathbb{P}^1$ we obtain a space of deformations of $X$ of
dimension $h^{1,2}$.

To find a space of deformations for $\hat{X}$ we perform a
simultaneous resolution of singularities in the family by blowing
up the families of components of singular fibers.

We have found a space of  dimension $h^{1,2}(\hat{X})$ describing
deformations of $\hat{X}$. It has to be an open set in the
Kuranishi space as the latter is smooth. By Corollary \ref{generic
F not I1,II } the generic element of this set has no fibers of
type $F\times I_0 $ for $F\not=I_1 \text{ or } II$.  This ends the
proof.

\end{proof}

From the above proof we can find a geometric description of
deformations of $\hat{X}$ looking only at the singular fibers of
$S_1$ and $S_2$ and their incidence conditions. Indeed each fiber
of type $F\times I_0$ can be deformed in the deformation family
into $\chi(F)$ fibers of type $I_1 \times I_0$. A singular fiber
of type $I_n \times F$ can be deformed into the sum of one
singular fiber of type $I_n \times I_m$ with the same number of
components as the original singular fiber and a fiber of type $I_1
\times I_0$. More precisely:
\begin{itemize}
 \item Fibers of type $IV\times I_n$ can deform into a fiber of
type $I_3\times I_n$ and a fiber of type $I_1 \times I_0$.
 \item Fibers of type $III\times I_n$ can deform into a fiber of
types $I_2\times I_n$ and a fiber of type $I_1 \times I_0$.
\end{itemize}
Moreover, fibers of type $I_n\times I_m$ have to preserve their
type during the deformation.
\begin{ex}
To ilustrate the result of the above theorem let us take two surfaces $S_1$ and $S_2$ given by Weierstrass equations as in the proof of above Lemma with:
\begin{itemize}
 \item $f=12(t^4-t^2+1)$
 \item $g=4(2t^6-3t^4-3t^2+2)$
 \item $f'=3(t-3)^2(t-1)^2$
 \item $g'=(t-3)^2(t-1)^2((t-2)^2+1)$
\end{itemize}
The discriminants of the surfaces are
\begin{itemize}
 \item $\Delta=2^4 3^6 t^4(t+1)^2(t-1)^2$, 
 \item $\Delta'=-108 (t-3)^4 (t-1)^4 (t-2)^2$
\end{itemize}
In this case $S_1$ has singular fibers of types $I_2$ $I_2$ $I_4$ $I_4$ in the points $-1,1,0,\infty$ and $S_2$ singular fibers of types $I_2$ $I_2$ $IV$ $IV$ in the points $2,\infty,3,1$. Their fiber product has one fiber of type $I_2 \times IV$ and one fiber of type $I_4\times I_2$. Remaining fibers are of types $F \times I_0$. As described in the the proof of above Lemma the conditions imposed on quadruples $F,F',G,G'$ describing the set $Q$ come from the requirement that the discriminants of both Weierstrass equations have two common zeroes of the following types. One of the common zeroes needs to be double for the frist discrminant and triple for the second and the second zero needs to be quadruple for the first discriminant and double for the second. These are 5+6-2 conditions on the space of discriminants hence also on $Q$. This gives $Q$ of dimension 15. Finally the deformation space is of dimension 10. 
\end{ex}

\begin{rem} We did not prove that a fiber $F\times I_n$, where $F$
is an unstable fiber always splits in the deformation space. Using
\cite{Perlist} we can prove that it is always possible to split a
fiber of types $II$ (resp. $III$, $IV$) into a sum of fibers $I_1$
and $I_1$ (resp. $I_2$ and $I_1$; $I_3$ and $I_1$). However it is
not clear that we can do it controlling the position of fibers.
\end{rem}

\begin{rem} In all above we did not need to assume that all fibers
of $S_1$ and $S_2$ are reduced. We can allow fibers of type
$F\times I_0$ for any $F$. Such fiber will also split in the
deformation space into fibers of type $I_1 \times I_0$.
\end{rem}
\subsection{The fiber product as an iterated double covering}
Here, for a given fiber product of two elliptic surfaces with
section we give another smooth model. This enables us to find an
argument to compute the deformation space of $\hat{X}$
independently of the characteristic of the base field. Let us
start with the following lemma concerning elliptic surfaces.

\begin{lem} \label{powierzchnia jako podwojna kwartyka }
Let $S$ be a rational elliptic surface with section. Assume
moreover that all fibers of $S$ are reduced. Then $S$ is
birational to the double covering $\hat{S}$ of $\mathbb{P}^2$
branched along a quartic curve $Q$ with only $A_k$ singularities.
Moreover, we can choose this birational equivalence to map fibers
of $S$ into the inverse images by the covering of lines passing
through a fixed point not lying on $Q$.
\end{lem}
\begin{proof} Let us choose two disjoint sections $s_1$ and $s_2$ of $S$.
To do this we use the fact that $S$ is the blow up of $\mathbb{P}^2$ in the
base locus of a system of cubics (see \cite{MirPer}). By the assumption that all
fibers are reduced we know that this system has at least two
non-infinitely near base points. Then each of the at least two -1
exceptional divisors of the blowing up of the base set is a
section and they are all pairwise disjoint.

Define an involution $i$ acting on fibers of $S$ such that
$i(s_1)=s_2$. This can be done in the following way, we define
$s_1$ to be the 0 section and take $i(x)=s_2-x$ on each fiber. The
quotient of $S$ by $i$ is a ruled surface with section $s$ (the
image of $s_1$ and $s_2$). After blowing down all components of
the fibers that are disjoint from $s$ we obtain a minimal ruled
surface $R$ with a $-1$ section and no section with self
intersection less than $-2$. But on the Hirzebruch surface $\mathbb{F}_2$ we do not have any $-1$
curves. Hence $R$ is the Hirzebruch surface $\mathbb{F}_1$ and $S$ is birational to the double
cover of this $\mathbb{F}_1$ branched in a four-section disjoint from the
section $s$. Blowing down the section $s$ we obtain the assertion
of the Lemma. Indeed all singularities of the quartic obtained in
this way are of type $A_k$ as resolving triple points of the
branch locus lead to non-reduced fibers.
\end{proof}
\begin{rem}
In fact the only elliptic surfaces with no two disjoint sections
are the two surfaces $X_{22}$ and $X_{211}$ from \cite{MirPer}
(the only rational elliptic surfaces with a fiber of type
$I_2^*$). Hence we can generalize this lemma to all rational
elliptic surface with section except the two above ones, we need
only to allow all ADE singularities.
\end{rem}
\begin{rem} The above lemma and remark lead us to an interesting
correspondence between elliptic surfaces with section (except
$X_{22}$ and $X_{211}$) and Gorenstein del Pezzo surfaces of
degree 2 (see Chapter 4).
\end{rem}

For clarity let us consider the case when $S_1$ and $S_2$ have
only semi-stable fibers. Let $\hat{S}_1$ and $\hat{S}_2$ be
surfaces defined in Lemma \ref{powierzchnia jako podwojna kwartyka
}. The fiber product $S_1 \times_{\mathbb{P}^1} S_2$ is birational
to the variety given in the product $\mathbb{P}(1,1,1,2)\times
\mathbb{P}(1,1,1,2)$ with coordinates $(a,b,c,d;\alpha, \beta,
\gamma, \delta)$ by the equations
$$ d^2=Q_1(a,b,c), \quad \delta^2=Q_2(\alpha,\beta,\gamma), \quad a=\alpha, \quad b=\beta.
$$
This is birational to the variety $Y$ constructed as follows. Take
the double cover of $\mathbb{P}^3$ branched over a quartic cone
$Q_1(x,y,z)=0$ and next the double cover of the obtained variety
branched over the inverse image of the cone $Q_2(x,y,t)=0$ by the
first covering. Thus the fiber product $S_1 \times_{\mathbb{P}^1}
S_2$ is birational to the iterated double cover $Y$. Observe that
$Y$ is a singular Calabi--Yau three fold with a smooth model
$\hat{Y}$ constructed in the following way. We blow up
$\mathbb{P}^3$ in such a way that the quartic cone $Q_1(x,y,z)=0$
is resolved, i.e. first we blow up the vertex of the cone then
consecutively all the double lines (also the infinitely near).
Next we take the first double covering branched over the proper
transform of the cone. After that we do the same with the inverse
image of the cone $Q_2(x,y,t)=0$ to the double cover. Here after
blowing up the fourfold points and the double curves we obtain two
new kind of singularities:
\begin{itemize}
\item Some nodes of the three fold lying outside the branch locus.
\item Some nodes of the branch locus.
\end{itemize}
Both induce nodes on the double covering that we resolve by taking
a small resolution.

The results of \cite{CV} give another method to compute the
dimension $h^{1,2}(\hat{X})=h^{1,2}(\hat{Y})$ of the space of
deformations of $\hat{X}$.

Let us use the method in our context. We prove the following
theorem

\begin{thm}\label{deformacje iterowanego} Let $Y$ be an iterated double cover of
$\mathbb{P}^3$ branched over two quartic cones $Q_1$ and $Q_2$
defining elliptic surfaces with semi-stable fibers. Assume
moreover that $Q_1$ ad $Q_2$ do not pass through each other
vertices. Then deformations of the small resolution of $Y$
correspond to those deformations of both cones that preserve the
tangency relations and the multiplicities of the intersection
points.
\end{thm}
\begin{ex}
Let us consider two reducible quartic cones in $\mathbb{P}^3$ (with coordinates $(x,y,z,t)$), given by equations $Q_1=\{x(x+z)((x-z)^2+t^2-25z^2)=0\}$ and $Q_2=\{y(y+z)((x-3z)^2 +t^2-25z^2)=0\}$. These two cones are tangent in the points $(1,3,1,5)$ and $(1,3,1,-5)$. They moreover intersect each other nontranversely in the point $(0,0,0,1)$. It is worth pointing out that the remaining points of intersection  
including the points of intersections of double lines on one cone with the other cone are treated as transversal, hence of multiplicity one. The above theorem says that deformations of a small resolution $\tilde{Y}$ of $\mathbb{P}^3$ branched over the two cones $Q_1$ and $Q_2$ correspond to such deformations of both cones that preserve tangency in one point and one quadruple point of intersection. More explicitly a small deformation of $\tilde{Y}$ is an iterated cover of $\mathbb{P}^3$ branched over two quartic cones tangent in one point and admitting double lines intersecting in one point. 
\end{ex}

\begin{cor} Any deformation of a fiber product of two elliptic
surfaces with section admitting only semi-stable fibers is a fiber
product of elliptic surfaces with section with corresponding
singular fibers consisting of the same number of components as in the
deformed variety.
\end{cor}
\begin{proof}
We have already proved that a fiber product of elliptic surfaces with section is a small resolution of an iterated double cover of $\mathbb{P}^3$ branched over two quartic cones. By Theorem \ref{deformacje iterowanego}, deformations of such varieties correspond to deformations of the cones. The deformed cones induce new fiber products. The only thing we need to observe is that a fiber of type $I_n \times I_m$ arise only in the case when the cones are either tangent to each other, or a double line of one cone is tangent to the other, or two double lines of the cones intersect. From the latter and Theorem \ref{deformacje iterowanego} it follows that the fibers of type $I_n \times I_m$ are preserved in deformations, which is the assertion of the corollary.
\end{proof}

\begin{proof}[Proof of theorem \ref{deformacje iterowanego}]
Let us start with the observation that the space of deformations
of a Calabi--Yau three fold which is a small resolution of a
singular Calabi--Yau three fold $Y$ is isomorphic to the space of
deformations of the variety which is a big resolution of $Y$. We
study deformations of the three fold obtained by taking the big
resolutions of the nodes.

We introduce the following notation:
\begin{itemize}
 \item Let $D_1$ and $D_2$ denote the quartic cones from the
 theorem.
 \item Let $\sigma_1:\tilde{\mathbb{P}}^3\ra \mathbb{P}^3$
be a sequence of blowups of $\mathbb{P}^3$ inducing a minimal
resolution of $D_1$. The proper transforms of $D_1$ and $D_2$ by
$\sigma_1$ will be denoted by $\tilde{D}_1$ and $\tilde{D}_2$. The
line bundles corresponding to $\frac{1}{2}\tilde{D}_1$ and
$\frac{1}{2}\tilde{D}_2$ will be denoted $\tilde{L}_1$ and
$\tilde{L}_2$.
 \item Let $\varphi_1:\tilde{Z}\ra\tilde{\mathbb{P}}^3$ be the double cover of
 $\tilde{\mathbb{P}}^3$ branched over $\tilde{D}_1$ and given by
 the corresponding line bundle $\tilde{L_1}$. The pullback of the
 bundle $\tilde{L}_2$ will be denoted by $\tilde{\mathcal{L}}_2$.
 \item Let $\sigma_2:\hat{Z}\ra\tilde{Z}$ be a sequence of blow ups of $\tilde{Z}$
 inducing a minimal resolution of $\varphi_1^{-1}(\tilde{D}_2)$.
 The proper transforms of $\varphi_1^{-1}(\tilde{D}_1)$ and
 $\varphi_1^{-1}(\tilde{D}_2)$ by $\sigma_2$ will be denoted by
 $\hat{D_1}$ and $\hat{D}_2$. The line bundle corresponding to $\frac{1}{2}\hat{D}_2$ will
 be denoted by $\hat{\mathcal{L}}_2$.
 \item Let $\varphi_2:\hat{Y}\ra\hat{Z}$ be the double cover of
 $\hat{Z}$ branched over $\hat{D_2}$.
\end{itemize}
We have the following diagram
$$
\begin{CD}
      Y    @<<<               \tilde{Y}     @ <<<         \hat{Y}     \\
 @   VV2:1V   @                     VV2:1V        @              V\varphi_{2}V2:1V        \\
     Z     @<<<               \tilde{Z}     @ <\sigma_2<<         \hat{Z}     \\
 @  VV2:1V    @                 V\varphi_{1}V2:1V            @.                        \\
\mathbb{P}^3@<<\sigma_1< \tilde{\mathbb{P}}^3 @.                        \\
\end{CD}
$$
We are interested in computing $h^1(\hat{Y},\Theta_{\hat{Y}})$. As
$\varphi_2$ is a finite morphism we have
$$h^1(\hat{Y},\Theta_{\hat{Y}})=h^1(\hat{Z},
\varphi_{2*}(\Theta_{\hat{Y}}))= h^1(\hat{Z},\Theta_{\hat{Z}}
\otimes \hat{\mathcal{L}}_2^{-1})+
h^1(\hat{Z},\Theta_{\hat{Z}}(\mathrm{log} \hat{D}_2 )).$$

The first component of the sum can be computed using the fact that
$\hat{Y}$ is a big resolution of a nodal Calabi--Yau three fold.
We have then
$$h^1(\hat{Z},\Theta_{\hat{Z}}
\otimes \hat{\mathcal{L}}_2^{-1})=
h^1(\tilde{Z},\Theta_{\tilde{Z}} \otimes
\tilde{\mathcal{L}}_2^{-1})+ \sum g(\hat{C}_i), $$
 where $g(\hat{C}_i)$ are the genera of the blown up curves. Each
 of these curves is a component of the pre-image of some double curve
 $C_i$ on the cone $D_2$. Hence their genera are equal to one if the curve
 $C_i$ was not tangent to the cone $D_1$ and zero otherwise.
Taking again the direct image, this time by $\varphi_1$ we obtain
by the projection formula

$$h^1(\tilde{Z},\Theta_{\tilde{Z}} \otimes
\tilde{\mathcal{L}}_2^{-1})=
h^1(\tilde{\mathbb{P}}^3,\Theta_{\tilde{\mathbb{P}}^3} \otimes
\tilde{L}_1^{-1}\otimes
\tilde{L}_2^{-1})+h^1(\Theta_{\tilde{\mathbb{P}}^3}(\mathrm{log}\tilde{D}_1)
\otimes \tilde{L}_2^{-1})=0$$

We need to prove the last equality.

\begin{prop} The following equalities hold.
\begin{itemize}
 \item $h^1(\tilde{\mathbb{P}}^3,\Theta_{\tilde{\mathbb{P}}^3}
\otimes \tilde{L}_1^{-1}\otimes \tilde{L}_2^{-1})=0,$
 \item $h^1(\Theta_{\tilde{\mathbb{P}}^3}(\mathrm{log}\tilde{D}_1)
\otimes \tilde{L}_2^{-1})=0.$
\end{itemize}
\end{prop}
\begin{proof} We follow the idea of \cite[section 5]{CV}. We have the exact sequence
\begin{equation}\label{exact seq logD1}
0\longrightarrow
\Theta_{\tilde{\mathbb{P}}^3}(\operatorname{log}(\tilde{D}_1))
\otimes \tilde{L}_2^{-1} \longrightarrow
\Theta_{\tilde{\mathbb{P}}^3} \otimes \tilde{L}_2^{-1}
\longrightarrow N_{\tilde{D}_1 | \tilde{\mathbb{P}}^3}\otimes
\tilde{L}_2^{-1} \longrightarrow 0.
\end{equation}
We prove now some lemmas.
\begin{lem} The following equality holds.
$$h^0(N_{\tilde{D}_1 | \tilde{\mathbb{P}}^3}\otimes
\tilde{L}_2^{-1})=0$$
\end{lem}
\begin{proof}[Proof of Lemma] The divisor $(\tilde{D}_1\otimes
\tilde{L}_2^{-1})|_{\tilde{D}_1}$ is not effective as its index of
intersection with the proper transforms of the rays of the cone
$D_1$ is negative.
\end{proof}

\begin{lem} The following equalities hold. $$h^1(\Theta_{\tilde{\mathbb{P}}^3} \otimes
\tilde{L}_2^{-1})=h^1(\tilde{\mathbb{P}}^3,\Theta_{\tilde{\mathbb{P}}^3}
\otimes \tilde{L}_1^{-1}\otimes \tilde{L}_2^{-1})=0$$
\end{lem}
\begin{proof}[Proof of Lemma] We use the Leray spectral sequences associated to
the blowings up. The number
$h^1(\tilde{\mathbb{P}}^3,\Theta_{\tilde{\mathbb{P}}^3} \otimes
\tilde{L}_1^{-1}\otimes \tilde{L}_2^{-1})$ is computed explicitly
in \cite[section 5]{CV} as $\tilde{\mathbb{P}}^3$ is a partial
step to resolving the octic $D_1+D_2$. To compute
$h^1(\Theta_{\tilde{\mathbb{P}}^3} \otimes \tilde{L}_2^{-1})$ we
use the same argument as in the cited paper. We proceed by
recursion. Let $\sigma: \tilde{P}\rightarrow P$ be a single blow
up from the sequence leading from $\mathbb{P}^3$ to
$\tilde{\mathbb{P}^3}$. As $D_2$ does not pass through any center
of the blowings up $\tilde{L}_2$ can be replaced in each step of
the recursion by $\sigma^*(L_2)$. We use the equality
$$h^1(\Theta_{\tilde{P}} \otimes \sigma^*(L_2)^{-1})=
h^1(\Omega_{\tilde{P}}^{2} \otimes K_P^{\vee}\otimes
\sigma^*(L_2)^{-1})=h^{2}(\Omega_{\tilde{P}}^{1} \otimes
K_P\otimes \sigma ^*(L_2))$$
 and the exact sequence
\begin{eqnarray*}
0\longrightarrow \sigma^*(\Omega_{P}^{1} \otimes K_P\otimes L_2)\otimes
\mathcal{O}_{\tilde{P}}(k E)\longrightarrow \Omega_{\tilde{P}}^{1}
\otimes K_P\otimes \sigma^*(L_2) \longrightarrow\\
\longrightarrow \Omega^1_{E/C}
\otimes \mathcal{O}_E (-k)\otimes \sigma^*(K_P\otimes
L_2)\longrightarrow 0,
\end{eqnarray*}
where $C$ is the blown up set, $k=\operatorname{codim}_P (C)-1$,
 and $E$ is the exceptional divisor. We have
$k>0$. From the projection formula we deduce
 $$\sigma_*(\Omega_{\tilde{P}}^{1} \otimes
K_P\otimes \sigma ^*(L_2))=\Omega_P^{1} \otimes K_P \otimes L_2
\otimes \sigma_*(\mathcal{O}_{\tilde{P}}(k E))=\Omega_P^{1}
\otimes K_P \otimes L_2.$$
 From the derived exact sequence we obtain $R^i \sigma_*(\Omega_{\tilde{P}}^{1} \otimes
K_P\otimes \sigma ^*(L_2))=0$ for $i\geq 1$, since $R^i
\sigma_*(\Omega_{E/C}^{1}(-k)=0$ for $i\geq 0$ and $R^i \sigma_*(
\mathcal{O}_{\tilde{P}}(k E))=0$ for $i\geq 1$. Using the Leray
spectral sequence and the recursion we get at the end
$h^1(\Theta_{\tilde{\mathbb{P}}^3} \otimes
\tilde{L}_2^{-1})=h^1(\Theta_{\mathbb{P}^3} \otimes L_2^{-1})=0.$
\end{proof}
The exact sequence \ref{exact seq logD1} together with the lemmas
proves the proposition.
\end{proof}

The second component $h^1(\hat{Z},\Theta_{\hat{Z}}(\mathrm{log}
\hat{D}_2 ))$ of the sum can be interpreted as the dimension of
the space of deformations of $\hat{D}_2$ contained in deformations
of $\hat{Z}$. It is thus the space of equisingular deformations of
$\tilde{D}_2$ in deformations of $\tilde{Z}$. Let us first find a
geometric interpretation of the space of all small deformations of
$\tilde{Z}$. To do this we compute as above
$$h^1(\tilde{Z}, \Theta_{\tilde{Z}})=h^1 (\tilde{\mathbb{P}}^3,
\Theta_{\tilde{\mathbb{P}}^3}\otimes \tilde{L}_1^{-1})+ h^1(
\Theta_{\tilde{\mathbb{P}}^3}(\mathrm{log} D_1)).$$
 Now, $h^1(\Theta_{\tilde{\mathbb{P}}^3}(\mathrm{log} D_1))$
 represents the space of equisingular deformations of $D_1$ in
 $\mathbb{P}^3$. These deformations are also cones as they are of
 degree 4 and admit a fourfold point. Moreover, they correspond to
 equisingular deformations of the base quartic of $D_1$ in
 $\mathbb{P}^2$. The dimension of the space of transverse deformations
 is given by
 $$h^1 (\tilde{\mathbb{P}}^3,\Theta_{\tilde{\mathbb{P}}^3}\otimes \tilde{L}_1^{-1})= \text{The
 number of blowings up of double curves}.$$
The above two spaces generate together the space of all quartic
cones in $\mathbb{P}^3$. Hence, every deformation of $\tilde{Z}$
is a resolution of the double cover $Z$ of $\mathbb{P}^3$ branched
over a quartic cone. Next, observe that a deformation of
$\tilde{D}_1$ in a deformation of $\tilde{Z}$ corresponds to a
deformation of its projection onto $Z$. The latter is a
deformation of a complete intersection of two quartics in
$\mathbb{P}(1,1,1,2)$ hence is also an intersection of this type.
Thus equisingular deformations of $\tilde{D}_1$ in deformations of
$\tilde{Z}$ are a subset of deformations of a pair of cones in
$\mathbb{P}^3$ with fixed vertices. It remains to observe that the
equisingularity condition imposes on the deformed cones the same
tangency and incidence conditions as in the original cones.
\end{proof}
\begin{rem} All above works also if we allow any reduced fibers.
The only difference is that during the resolution of the branch
locus of the second covering we obtain singularities that are not
nodes. These are some double and triple points admitting small
crepant resolutions. As these are only double and triple points
they do not affect the space of transversal deformations. The rest
remains unchanged.
\end{rem}
\begin{rem} We can use the same argument to compute deformations of
the fiber product of general elliptic surfaces with section.
Indeed, most non reduced fibers correspond to triple curves on the
cones. We treat the triple curves in the following way. We blow
them up and add the exceptional divisor to the proper transform of
the branch locus on the blowing up to obtain an even divisor. The
rest remains unchanged.
\end{rem}

\section{Kummer fibrations}
In this section we study a fiberwise Kummer construction for a
product of rational elliptic surfaces with section.

Let $X$ be the fiber product (possibly singular) of two rational
elliptic surfaces $S_1$ and $S_2$ with section. Let
$a=(a_1,a_2):\mathbb{P}^1\longrightarrow X $ be any section on
$X$. Then $a_1$ and $a_2$ are sections of $S_1$ and $S_2$
respectively. Let us consider the involution $i: X\longrightarrow
X$ such that on each smooth fiber it is of the form
$i(x_1,x_2)=(a_1-x_1, a_2 -x_2)$. Observe that $i$ is well defined
as we have a group structure on each fiber. Let $Y$ be the
quotient of $X$ by the involution $i$.
\begin{prop} \label{kummera jest CY} In the above setting $Y$ admits a resolution being a
Calabi-Yau three fold (not necessarily projective).
\end{prop}

Before the proof of the proposition let us consider more closely
the involution $i$. Let $i_1$ and $i_2$ denote the involutions on
$S_1$ and $S_2$ such that $i=(i_1,i_2)$. We study the possible
actions of these involutions on a singular fiber $F$ of $S_1$ or
$S_2$ respectively.

\begin{lem} Let $S$ be a rational elliptic surface with chosen 0
section and such that $S$ admits only reduced fibers. Let $j$ be
an involution of the form $x \mapsto b-x$, where $b$ is a section
of $S$. Let $F$ be a singular fiber of $S$. Then exactly one of
the following possibilities occurs.
\begin{enumerate}
 \item The fiber $F$ is of type $I_1$ and $j$ acts on $F$ by
symmetry.
 \item The fiber $F$ is of type $I_{2k}$. Then we have one of the
 following cases.
\begin{enumerate}
\item The involution $j$ has two fixed points (these are then two
opposite nodes) and interchanges pairs of components of $F$.

\item The involution $j$ acts on two opposite components of $F$ by
symmetry, interchanging respective remaining components.
\end{enumerate}
 \item The fiber $F$ is of type $I_{2k+1}$ for $k\geq 1$. Then $j$
acts on one of the components of $F$ by symmetry with fixed points
outside singularities of $F$ and interchanges the respective pairs
of remaining components.
 \item The fiber $F$ is of type $II$. Then $j$ fixes the
 singular point of $F$ and admit one more fixed point.
 \item The fiber $F$ is of type $III$. Then we have two
 possibilities:
\begin{enumerate}
 \item The involution $j$ has only one fixed point and interchanges
the two components of $F$.
 \item The involution $j$ acts on both components fixing their
 intersection point and one more point on each fiber.
\end{enumerate}
 \item The fiber $F$ is of type $IV$. Then the involution $j$ fixes
 the triple point, interchanges two of the components
 of the fiber and acts on the third component, admitting an
 additional fixed point.

\end{enumerate}
\end{lem}

\begin{proof}[Proof of lemma]
Following \cite{Shaf2} for each fiber $F$ of $S$ there is an
induced group structure on $F^{\sharp}$, the set of smooth points
of the fiber. The involution $j$ is of the form $x\mapsto b-x$ on
each fiber including the singular ones. Hence $j$ acts also in
this way on the group $F^{\sharp}/F^0$ of connected components of
$F^{\sharp}$. The latter group for a reduced fiber is isomorphic
to $\mathbb{Z}/n\mathbb{Z}$. On this group we have three
possibilities for an involution of the form $[x]\mapsto [b]-[x]$.
\begin{itemize}
 \item The number $n$ is odd. Then the involution has one fixed
point (i.e. $j$ acts on exactly one of the component of the
singular fiber)
 \item The number $n$ is even but $[b]$ is odd. Then the involution
 has no fixed points (i.e. $j$ interchanges pairs of components of the singular fiber).
 \item The numbers $n$ and $[b]$ are both even. Then the involution
 has two fixed points (i.e. $j$ acts on exactly two of the component of the
singular fiber).
\end{itemize}
Moreover, the structure group on $F^0$, the component of the zero
section, is multiplicative for $F=I_n$ and additive for the
remaining fibers. This means that an involution fixing a component
of the fiber $I_n$ admits two fixed points on the smooth part of
this fiber. In the same way an involution fixing a component of a
fiber different from $I_n$ admits one fixed point on the smooth
part of this fiber. Together we obtain all possibilities described
in the assertion of the lemma.
\end{proof}

The above lemma implies immediately the following

\begin{cor}\label{cor o postaci lokalnej inwolucji}
Locally around a fixed point $P$ of $j$ we have the following
possibilities.
\begin{enumerate}
\item The fixed point $P$ is a smooth point of the fiber on which
it is lying.

\item The fiber admits a node in $P$ and the local analytic
components of this node are interchanged by the involution $j$.

\item The fiber is of type $III$ singular at $P$ and the components
of this fiber are interchanged by the involution $j$.

\item The fiber is of type $III$. The point of intersection is
fixed but the involution acts on both components separately.

\item The fiber is of type $IV$ and all three of its components pass
through $P$. Two of them are interchanged by the involution.
\end{enumerate}
\end{cor}

\begin{cor} \label{cor o rownaniach lok osobliw wlokn}
We can perform a local analytic change of coordinates around $P$
in such a way that the fibration is preserved and the surface $S$
is given in $\mathbb{C}^3$ with coordinates $(x,t,u)$ by the
equation:
\begin{enumerate}
 \item $u^2= x $ and the fibration is given by $t$,
 \item $u^2= x^2-t$ and the fibration is given by $t$,
 \item $u^2= x^4-t$ and the fibration is given by $t$,
 \item $u^2=x-t$ and the fibration is given by $xt$,
 \item $u^2=x^2 -t$ and the fibration is given by $xt$.
\end{enumerate}
Where the involution is given by $u\mapsto -u$.
\end{cor}
\begin{proof}
We need only to observe that each of the cases from the corollary
\ref{cor o postaci lokalnej inwolucji} is represented locally by
one of these equations.
\end{proof}

Let us come back now to the proposition.
\begin{proof}[Proof of Proposition \ref{kummera jest CY}]
By the adjunction formula for a double covering, as $Y$ is a double cover of a Calabi--Yau three fold branched over a set of codimension 2, the canonical divisor $K_Y$ is trivial.
We need only to find a crepant resolution of the singularities.
Outside fixed points of the involution the singularities are
locally isomorphic to the corresponding singularities of $X$,
hence admit a crepant resolution. Moreover, from the Corollary
\ref{cor o rownaniach lok osobliw wlokn} around the fixed points
of the involution $i$ we have the following possibilities for the
local equation of the quotient variety.
\begin{enumerate}

 \item[(a)] $u^2=xy$
 \item[(b)] $u^2=(x^2-t)y$
 \item[(c)] $u^2=(x^4-t)y$
 \item[(d)] $u^2=(x-t)y$
 \item[(e)] $u^2=(x^2-t)(y^2-t)$
 \item[(f)] $u^2=(x^4-t)(y^2-t)$
 \item[(g)] $u^2=(x-t)(y^2-xt)$
 \item[(h)] $u^2=(x^4-t)(y^4-t)$
 \item[(i)] $u^2=(x^2-t)(y^2-xt)$
 \item[(j)] $u^2=(x-t)(y^4-xt)$
 \item[(k)] $u^2=(x-t)(y-z)$, $xt=yz$

\end{enumerate}
The first step in the resolution of the quotient variety is the
blowing up of the branch curve. Observe that the first four
possibilities are transversal $A_1$ singularities. It means that
the above blow up gives their crepant resolution. In all remaining
cases except case (k) after blowing up the singular curve we get
isolated singularities with small resolution.
 \begin{itemize}
\item[(e)] We get two nodes.

\item[(f)] We get two singularities of type $A_3$.

\item[(g)] We get two singularities of type $D_4$.

\item[(h)] We get two singularities given in suitable local
coordinate systems by the equation $x^4-y^4=u^2-t^2$.

\item[(i)] We get a singularity of type $A_5$.

\item[(j)] We get a singularity given in a suitable local
coordinate system by the equation $x^4-y^4=u^2-t^2$.
\end{itemize}

In the last case after blowing up the singular curve we obtain a
new singular curve that we next blow up again. As the blown up
curves were double curves the obtained resolution is also crepant.
\end{proof}
 From now on we assume moreover that $S_1$ and $S_2$ admit only semi-stable singular fibers.
This means that in the above list of local equations of singularities of the Kummer fibration we allow only types (a),(b) and (e).
\begin{prop} If $X$ has only semi-stable fibers and no fibers of type $I_1 \times I_n$,
then the three fold $Y$ admits a projective Calabi--Yau
resolution.
\end{prop}
\begin{proof}Observe first that outside the fixed locus $X^i$ of the
involution the proposition is trivial as the double cover
$X\setminus X^i\longrightarrow Y\setminus X^i$ is \'etale.

It remains to see that blowing up separately the strict transforms
of all components of all fibers we can resolve the nodes that
appeared on the partial resolution (after blowing up the double
curves) of $Y$. Here the only case we need to check more closely
is what happens over a singular point of the fifth type from the
proof of Proposition \ref{kummera jest CY}. After blowing up the
branch curve (locally it is defined by the equations
$x^2-t=y^2-t=u=0$) the two components of the fiber $t=0$ (given
locally by equations $t=0$, $u=\pm xy$) remain smooth and pass
through the obtained nodes. The fact that these are indeed two
global components follows from the assumption that fibers are
semi-stable and no fiber is of type $I_1\times I_n$.
\end{proof}
\subsection{Hodge numbers}
We compute Hodge numbers of the constructed projective varieties.
Let us start by computing their Euler characteristics. Let us
denote the branch curve of the involution $i$ by $C$. From the
above construction we have
$$ \chi(\hat{Y})= \frac{\chi(X)-\chi(C)}{2} + 2\chi(C) + 2o + b ,$$
where $o$ denotes the number of singular points of $X$ fixed by
the involution and b denotes the half of the number of the
remaining nodes of $X$. The number $\chi(C)$ can be computed from
the number of fixed points on each fiber.
\begin{lem} \label{chi c} We have the following equality:
$$\chi(C)=16(2-\sharp(A_1 \cup A_2)) + \sum_{t\in A_1 \cup A_2} a_t a'_t,$$
where $a_t$ and $a'_t$ are the number of fixed points of the
involutions $i_1$ and $i_2$ on fibers lying over $t$ of the
respective surfaces $S_1$ and $S_2$.
\end{lem}
\begin{proof} The proof is straightforward as $C$ admits a
projection onto $\mathbb{P}^1$ with generic fiber consisting of 16
points. Special fibers appear only on singular fibers of the
fibration.
\end{proof}
\begin{rem}
Depending on the type of fibers the numbers $a_t$ and $a'_t$ take
only one of the three values \{2,3,4\}. Moreover, $a_t=3$ if and
only if the fiber $F_t$ of $S_1$ is of type $I_n$ for n odd.
\end{rem}
\begin{rem} The Euler characteristic of $\hat{Y}$ depends on the
involution $i$ i.e. for some Calabi--Yau fiber products we can
produce Kummer three folds with different Hodge numbers.
\end{rem}
The remark is illustrated by the following example.
\begin{ex}
Let $S_1$ and $S_2$ be double covers of $\mathbb{P}^2$ branched
over generic quartic curves. Let their elliptic fibrations be
chosen in such a way that $S_1$ has exactly one fiber of type
$I_2$ (one of the lines defining a fiber is double tangent to the
quartic), all remaining singular fibers of $S_1$ are of type
$I_1$, and all singular fibers of $S_2$ are of type $I_1$. Let the
0 sections of both surfaces be given by the pre-image of the
points defining the fibrations. We denote by $i$ the involution
defined by the zero section of the product and by $j=(j_1,j_2)$ the
involution induced by the coverings defining each surface. The involution $j=(j_1,j_2)$ is also an involution of type $(a_1-x_1,a_2-x_2)$ because both $j_1$ and $j_2$ have four fixed points on the generic fiber of $S_1$ and $S_2$. The sections $a_1$ (resp. $a_2$) is the section for which the preimages of the branch quartic defining $S_1$ (resp. $S_2$) by the double covering represent the four-section given by the equation $2x_1=a_1$ (resp. $2x_2=a_2$) in the group structure of each smooth fiber of $S_1$.     As we
have chosen everything to be generic, the product $S_1
\times_{\mathbb{P}^1}S_2$ is already smooth, hence $o=b=0$. Both
involutions on each fiber of the fibration different from the
fiber $I_2 \times I_0$ have the same number of fixed points.
However, on this fiber they differ ($i$ has 4 fixed points but
$j$ only two) giving different $\chi(C)$ and consequently
different $\chi(\hat{Y})$. More precisely $\chi(\hat{Y^i})= -74$
and $\chi(\hat{Y^j})=-96$.

\end{ex}

Let us now concentrate on computing the deformations of $\hat{Y}$.
To do so we will need some more results concerning elliptic
surfaces with involution.

\begin{lem} Let $S$ be a rational elliptic surface with a chosen 0 section. Let $b$
be any section of $S$. Let $j$ be the involution given in each
smooth fiber by $x\mapsto b-x$. Then at least one of the following two possibilities hold:
\begin{itemize}
\item The surface $S$ is birational to the double cover of
$\mathbb{P}^2$ branched in a quartic curve with ADE singularities.
This birational equivalence can be chosen in such a way that the
fibration of $S$ is given by proper transforms of lines passing
through a point on $\mathbb{P}^2$ and that the involution $j$
corresponds to the involution of the covering.

\item The surface is birational to the double cover of a quadric
cone in $\mathbb{P}^3$ (i.e. $\mathbb{P}(1,1,2)$) branched in the
vertex and the intersection of this cone with a smooth cubic not
passing through the vertex. The birational equivalence can be
taken in such a way that the fibration is given by the proper
transforms of the rays of the cone and the involution $j$
corresponds to the involution of the covering.
\end{itemize}
\end{lem}
\begin{proof}
Observe that taking the quotient of the surface $S$ by the
involution $j$ we obtain a ruled surface $R$ over $\mathbb{P}^1$. This
ruled surface is either $\mathbb{P}^1\times \mathbb{P}^1$ or it
admits a section $s$ which is an exceptional curve. 

Let us first consider the second possibility. Choose $s$ and
blow down all components of the fibers disjoint from it. In this
way we obtain a minimal ruled surface with section $s$. Note that each component of
the pullback of $s$ by the covering is a -1 curve and a section on the elliptic
surface. Now taking into account the position of $s$ with respect to the branch divisor we can compute its self-intersection number. We have three possibilities:
\begin{itemize}
 \item The section $s$ is contained in the branch locus of the
 quotient by the involution. Then the pullback of $s$ by the covering is a double -1 curve. The self-intersection number of $s$ is then equal to -2. Hence $s$ is the -2 curve on the ruled surface $\mathbb{F}_2$.
 \item The section $s$ is disjoint from the branch locus of the
 quotient by the involution. Then its preimage by the covering has two disjoint components with self-intersection -1. The self-intersection number of $s$ is then equal to -1. Hence $s$ is the -1 curve on the ruled surface
 $\mathbb{F}_1$.
 \item The section $s$ cuts the branch locus of the quotient by the
 involution. Then the preimage of $s$ by the covering is either irreducibe with even self-intersection or has two intersecting components. In both cases $s$ has to have non-negative self-intersection hence cannot be an exceptional curve. 
\end{itemize}
After blowing down $s$ the assertion follows.

The remaining possibility is that the ruled surface $R$ is $\mathbb{P}^1\times\mathbb{P}^1$. In this case we can perform an elementary transformation with center at a tangency point of some fiber with the branch locus on $\mathbb{P}^1\times\mathbb{P}^1$ to obtain $\mathbb{F}_1$. We then need only to blow down the exceptional section of $\mathbb{F}_1$ to obtain (a) from the assertion.
\end{proof}
\begin{rem} \label{prawie zawsze stozek} The above lemma is very similar to Lemma \ref{powierzchnia jako podwojna kwartyka }. The only difference is that we take care not only of the surface and the fibration but also of the chosen involution. This additional requirement forces us to consider one more possibility.
However in the above proof if at least one fiber of the
elliptic surface has more than two components, then we have a choice of the section $s$ and can choose it not to be contained in the ramification divisor. This means that the only cases where the elliptic surface with a chosen involution is not birational to the double cover of $\mathbb{P}^2$ branched over a quartic curve with its natural involution is when all fibers of the surface are irreducible and $b$ is the 0-section. 
\end{rem}

 Now we can describe the Kummer fibration
$\hat{Y}$.
\begin{cor} \label{rownania Kummera}
The three fold $\hat{Y}$ is birational to one of the following:
\begin{itemize}
\item[(a)] The resolution of the double cover of $\mathbb{P}^3$
branched over the sum $D$ of two possibly reducible quartic cones.
 \item[(b)] The resolution of the double cover of
$\mathbb{P}(1,1,1,2)$ branched over the sum of two possibly
reducible weighted cones of degrees 4 and 6.
 \item[(c)]The resolution of the double cover of
$\mathbb{P}(1,1,2,2)$ branched over the sum of two possibly
reducible weighted cones of degrees 6 and 6.
\end{itemize}
Moreover, the Kummer fibration is given by proper transforms of
planes passing through the vertices of both cones.
\end{cor}

\begin{rem} Using Lemma \ref{powierzchnia jako podwojna kwartyka }
for every fiber product $X$ we can always construct an involution
$i$ on $X$ such that the corresponding Kummer fibration $\hat{Y}$
satisfies case (a).
\end{rem}

With regard to the above remark we consider only the case (a) of
Corollary \ref{rownania Kummera}. We comment later also the
remaining cases, however we omit details.

Let us assume that $\hat{Y}$ is projective and that the fixed
points set of the involution $i$ does not contain the zero
section. In this case $\hat{Y}$ is birational to the resolution of
the double covering of $\mathbb{P}^3$ branched along an octic
surface $D$, which is the sum of two quartic cones $Q_1$ and $Q_2$. From this point of view the natural smooth model to study
is described by the following construction. We blow up
$\mathbb{P}^3$ consecutively in all the fourfold points and the
double curves of the branch locus. Next, we take a small
resolution of all obtained nodes. This smooth model will be denoted by $\tilde{Y}$.

\begin{rem}
We have two different natural constructions of smooth models of introduced Kummer fibrations. The first one was described in the proof of Proposition \ref{kummera jest CY} the second is described above. In the generic case (e.g. the fiber product of two elliptic surfaces with singular fibers only of type $I_1$ and no common singular fibers) these two resolutions differ by the composition of two flops. More precisely in this case the first model arises from the second by blowing up each of the two components of the proper transform of the line passing through the vertices of the cones and contracting the second rulings of the obtained in this way quadrics.

In general the situation is more complicated. The described resolutions differ more when the elliptic surfaces admit common singular fibers which are fixed by the involution, i.e. when the two cones $Q_1$ and $Q_2$ are singular along some lines that intersect each other. 
\end{rem}
\begin{ex} To illustrate the second part of the above remark let us consider the following two elliptic surfaces. Let $S_1$ (resp. $S_2$) be the resolution of the double covering of $\mathbb{P}^2$ with coordinates $(x,z,t)$ (resp. $(y,z,t)$)  blown up in the point $(1,0,0)$, branched over the proper transforms of the nodal quartic $q_1=\{x(x+z)(x-t-z)(x+t-z)=0\}$ (resp. $q_2=\{y(y+z)(y-t-7z)(y+t+5z)=0\}$). Let their fibrations be given by the projection onto $(z,t)$. Observe that in this setting $S_1$ (resp. $S_2$) has singular fibers only of type $I_2$ in the points $\{(1,-2),(1,-1),(1,0), (1,1),(1,2),(0,1)\}$ (resp. $\{(1,-8),(1,-7),(1,-6),$ $(1,-5), (1,-4), (0,1)\}$). This means that teir fiber product admits one fiber of type $I_2 \times I_2$ in the point $(0,1)$ and all remaining fibers are of type $F\times I_0$.
The product is equipped with a natural involution coming from the involution of the covering on each surface. The Kummer fibration corresponding to this involuton is birational to the double covering of  $\mathbb{P}^3$ with coordinates $(x,y,z,t)$ branched over the sum of two quartic cones $Q_1=\{x(x+z)(x-t-z)(x+t-z)=0\}$ and $Q_2=\{y(y+z)(y-t-7z)(y+t+5z)=0\}$. Let us denote the vertices of the cones by $V_1=(0,1,0,0)$, $V_2=(1,0,0,0)$ and the remaining fourfold point of the octic surface $D=Q_1 \cup Q_2$ by $P=(0,0,0,1)$. The variety $\tilde{Y}$ introduced in the discussion above is obtained by blowing up $\mathbb{P}^3$ in the points $P$, $V_1$ and $V_2$, next blowing up the proper transforms of the lines $PV_1$ and $PV_2$ and the remainng doble lines of the cones, then blowing up the proper transform of the intersection $Q_1\cap Q_2$, and finally taking the double covering of the obtained variety branched over the proper tranform of $D$.  By straightforward computation in local coordinates we can prove that the resolution $\hat{Y}$ described in the proof of Theorem \ref{kummera jest CY} can be obtained from the resolution $\tilde{Y}$ by a sequence of four disjoint flops. Beside flopping the two components of the preimage by the covering of the proper transform of the line $V_1V_2$ we need to flop the two components of the preimage by the covering of the proper transform of the line of intersection of the exceptional divisor lying over $P$ with the proper transform of the plane $PV_1V_2$.
\end{ex}

Let us remind that the hodge numbers of Calabi--Yau threefolds are birational invariants. Hence to compute the deformations of
the manifold $\hat{Y}$ we can use results of \cite{CV} applied to $\tilde{Y}$.
The dimension of the space of transversal deformations of
$\tilde{Y}$ is equal to the geometric genus of the sum of all double
curves of the branch locus. As the double curves of each cone are
rational curves the only curve we need to deal with is the curve
of intersection of the cones denoted by $C_I$.
\begin{rem}
\label{transwersalne kummera} The curve $C_I$ is birational to the
branch curve $C$ of the involution on the fiber product. We
already know how to compute the Euler characteristic of the curve
$C$ (see Lemma \ref{chi c} ). Moreover, we know that $C$ is smooth
outside the nodes in the points of the fifth type from the proof
of Proposition \ref{kummera jest CY}. That is we know the Euler
characteristic of the resolution of $C$.
\end{rem}
This allows us to compute the geometric genus of $C_I$ provided we
also know the number of its components. Unfortunately this last
number cannot be deduced from the types of fibers of the fiber
product and the type of involutions induced on each of these
fibers. This is shown by the following example.
\begin{ex} \label{przyklad o ilosci skladowych krzywej rozgalezienia}
Let $S_1$ be the double cover of $\mathbb{P}^2$ branched over the
sum of a line and a nodal cubic cutting transversely. Let $S'_1$
be the double cover of $\mathbb{P}^2$ branched over a quartic with
three nodes. Let $S_2$ be the double cover of $\mathbb{P}^2$
branched over the sum of four generic lines. Let the fibration in
each case be given by preimages of lines passing through a generic
point of the plane. Then $S_1 \times_{\mathbb{P}^1} S_2$ has the
same types of fibers as $S'_1 \times_{\mathbb{P}^1} S_2$ and the
same type of involution (the type is given by the number of fixed
points) on each fiber. However the number of components of the
branch locus on the two three folds differ. In consequence the
dimensions of transversal deformations also differ. We will see
later that the space of equisingular deformations are equal in
both cases. Hence it is not enough to know the types of fibers of
the fibration and the involution on each fiber to compute the
Hodge numbers of a Kummer surface.
\end{ex}
\begin{rem} We can find a bound to the number of components by
studying the components of both cones. A lower bound is given by
the product of the number of components of both cones. Moreover,
with few exception this lower bound is in fact the actual number.
\end{rem}

By the methods of \cite{CV} it remains to compute the equisingular deformation of the surface $D$.
As we have assumed, $D$ is the sum of two quartic cones $Q_1$ and $Q_2$.
\begin{prop} \label{equisingularne kummera} 
The equisingular deformations of the branch locus $D$
correspond to a subset of deformations of the fiber product $\hat{X}$
consisting of those deformations that induce by Theorem \ref{deformacje iterowanego} equisingular
deformations of each of the quartic cones $Q_1$ and $Q_2$.
\end{prop}
\begin{proof}
We need first to prove that an equisingular deformation of $D$
induces a deformation of the fiber product $\hat{X}$.

To do this we observe that an equisingular deformation of the sum
of two quartic cones in $\mathbb{P}^3$ is also a sum of two
quartic cones. This follows from the fact that the curve of
intersection of the two cones is a complete intersection (4,4) and
hence has to be preserved in the deformation as a double curve
which is also a complete intersection of type (4,4). As we know
that the fourfold points are preserved the claim follows.

Next, we prove that the $I_n\times I_m$ fibers are preserved in
the equisingular deformations of $D$. We proceed as follows,
observe that a fiber of this type arises in one of the following
ways.
\begin{itemize}
\item The plane corresponding to the fiber is tangent to both
cones.
 \item The plane corresponding to the fiber is tangent to one of
the cones and contains some singular lines (one or two) of the
second cone.
 \item The plane corresponding to the fiber contains some singular
 lines of both cones.
\end{itemize}
Let us now observe that a double line is deformed into double
lines and a fourfold point has to be deformed into a fourfold
point. Thus the configuration of double lines is deformed with the
same incidence relations. Observe moreover that the intersection
points of lines tangent to both cones are also preserved in the
deformation as they induce two nodes on the blow up of
$\mathbb{P}^3$ in the curve of intersection of the cones. The same
concerns the intersection points of a double line with a tangency
line. This proves the claim as all fibers of type $I_n \times I_m$
are defined by the incidence relations between double and tangency
lines of the cones.

Together this gives an inclusion of the set of all equisingular
deformations of $D$ to the set of equisingular deformations of the
fiber product. The inclusion in the opposite direction follows
from the fact that we can perform a simultaneous resolution of
singularities in the family of all equisingular deformations of
the cones.

\end{proof}
\begin{rem}
Although general formulas describing the deformation space of a
studied Kummer fibrations are rather difficult to write, all above
computations are very easy for explicit examples.
\end{rem}
\begin{rem}
The above reasoning as it stands works only for Kummer fibrations
satisfying the case (a) from Corollary \ref{rownania Kummera}.
However the method of \cite{CV} should also work in the remaining
cases provided we prove that the surrounding weighted projective
spaces admit rigid resolutions. This is not hard as we have a
description of such spaces as a hypersurface and a complete
intersection respectively.
\end{rem}
\begin{rem} We have one more approach to deal with the cases (b)
and (c) from Corollary \ref{rownania Kummera}. We can use Remark
\ref{prawie zawsze stozek} to translate almost all Kummer
fibrations (except the products of a surfaces with only $I_1$
singularities) into resolutions of double covers of $\mathbb{P}^3$
branched over a sum of quartic cones. We need only to allow the
cones to pass through each other vertices. To this more general
picture we can also use results of \cite{CV}, the only difference
is that we have fivefold points, which are resolved in a more
complicated but also allowed way.
\end{rem}
\begin{rem} As earlier we can generalize the above results to non
projective Kummer fibrations.
\end{rem}

We know how to compute the deformations and the Euler
characteristic of the Kummer fibration $\hat{Y}$, hence we can
easily compute the rank of its Picard group.
\begin{ex} Let us take the most general fiber product. That is a
smooth variety $X$ which is a product of two elliptic surfaces
with only $I_1$ fibers (the singular fibers of the product are
only of type $I_1\times I_0$). The suitable Kummer fibration is
then the resolution $\hat{Y}$ of the double cover of
$\mathbb{P}^3$ branched over the sum of two generic quartic cones.
We have the following table of invariants.
\[
\begin{tabular}{c|c|c|c}
         &$\chi$ & $h^{1,2}$& $h^{1,1}$\\
\hline
      $X$& 0    & 19     & 19\\
$\hat{Y}$& -96  & 52& 4
\end{tabular}
\]
\end{ex}

\vskip10pt
\begin{minipage}{7cm} Grzegorz Kapustka \\
Jagiellonian University\\
ul. Reymonta 4\\
30-059 Krak\'ow\\
Poland\\
email:\\ Grzegorz.Kapustka@im.uj.edu.pl
\end{minipage}
\begin{minipage}{7cm}
 Micha\l\ Kapustka\\
 Jagiellonian University\\
 ul. Reymonta 4\\
 30-059  Krak\'ow\\
 Poland\\
 email:\\
 Michal.Kapustka@im.uj.edu.pl
\end{minipage}

\end{document}